\documentclass[12pt,twoside,reqno]{amsart}

\newtheorem{thm}{Theorem}[section]
\newtheorem{lemma}{Lemma}[section]
\newtheorem{remark}{Remark}[section]
\newtheorem{proposition}{Proposition}[section]

\def \R{{\Bbb R}}

\numberwithin{equation}{section}

\begin{document}

\title[generalized KdV equation]
{
	 Regular global solutions for a generalized KdV equation posed on  a half-line
}

\author[
N.~A. Larkin]
{
	N.~A. Larkin \\
	Departamento de Matem\'atica,\\
	Universidade Estadual de Maringá\\
	87020-900, Maringá, Parana, Brazil}
\


\thanks
{email:nlarkine@uem.br}

\begin{abstract}
An initial-boundary value problem for a generalized  KdV equation posed on   a half-line is considered.
Existence and uniqueness  of global regular solutions for arbitrary smooth initial data are established.
\end{abstract}

\maketitle

\section{Introduction}\label{introduction}

We are concerned with an initial-boundary value problem (IBVP) posed
on  the right half-line
$ x>0$  for the generalized KdV equation
\begin{equation}
Lu\equiv u_t+u^ku_x +u_{xxx}=0.\label{kdv}
\end{equation}

Equation (1.1) is a typical example of so-called dispersive equations attracting considerable attention
of both pure and applied mathematicians. The KdV
equation, $k=1,$  is  more studied in this context.
The theory of the initial-value problem
(IVP henceforth)
for (1.1) when $k=1$ is considerably advanced today
\cite{tao,kato,ponce2,saut2}.

 Although dispersive equations were deduced for the whole real line, necessity to calculate numerically the Cauchy problem approximating the real line by finite intervals implies to study initial-boundary value problems posed on bounded and unbounded intervals \cite{familark,larluc1,larluc2}.
What concerns (1.1) with $k>1, \;l=1$, called generalized KdV equations, the Cauchy problem  was studied in \cite{ martel,merle} and later in \cite{farah,Fonseka1,Fonseka2,ponce2}, where it has been established  that for $k=4$ (the critical case)  the problem is well-posed for small  initial data, whereas for arbitrary initial data solutions may blow-up in a finite time. The generalized Korteweg-de Vries equation was  studied for understanding the interaction between the dispersive term and the nonlinearity in the context of the theory of nonlinear dispersive evolution equations \cite{ponce2}. In \cite{lipaz}, the initial-boundary value problem for the generalized KdV equation with an internal damping posed on a bounded interval was studied in the critical case; exponential decay of weak solutions for small initial data has been established. \\
 Recently, due to physics and numerics needs, publications on initial-boundary value
problems in both bounded and unbounded domains for dispersive equations  have been appeared
\cite{larluc1,larluc2,saut2}. In
particular, it has been discovered that the KdV equation posed on a
bounded interval possesses an implicit internal dissipation. This allowed
to prove the exponential decay rate of small solutions for
\eqref{kdv} with $k=1$  posed on bounded intervals without adding any
artificial damping term \cite{larluc2}. Similar results were proved
for a wide class of dispersive equations of any odd order with one
space variable \cite{familark,larluc1,larluc2}.

 Our work has
been inspired by \cite{lipaz} where the critical KdV equation with internal damping
posed on a bounded interval   was considered and exponential decay of weak solutions has been established.

The main goal of our work is to prove  for arbitrary smooth initial data the existence and uniqueness
of global-in-time regular solutions for (1.1) posed  on the right half-line.
The paper is outlined as follows: Section I is the Introduction. Section 2 contains  formulation
of the problem and auxiliaries. In Section \ref{existence}, regularization by an initial--boundary value problem for the Kawahara equation
was used to prove the existence and exponential decay of regular global solutions without smallness conditions for initial data .

\section{Problem and preliminaries}\label{problem}

Let $T>0;$  and $\R^+= \{x\in \R^1,\; x>0\},\;\;Q_T=(0,T)\times \R^+.$
We use the usual notations of Sobolev spaces $W^{k,p}$, $L^p$ and $H^k$   and the following notations for the norms \cite{Adams}:

$$\| f \|_{L^p(\R^+)}^p = \int_{\R^+} | f  |^p\, dx,\;\;
\| f \|_{W^{k,p}(\R^+)} = \sum_{0 \leq | \alpha| \leq k} \|D^\alpha f \|_{L^p(\R^+)},\;p\in(1,+\infty).$$
$$H^k(\R^+)=W^{k,2}(\R^+);\;\;\|f\|_{L^{\infty}(\R^+)}=esssup_{\R^+}|f(x,t)|.$$

\vspace{8pt}

Consider the following IBVP:
\begin{align}
&Lu\equiv u_t+u^k u_x+u_{xxx}=0,\;\; \text{in}\;\; Q_T;\\
&u(0,t) =0,\; t>0;\\
&u(x,0)=u_0(x),\ \ x\in\R^+,
\end{align}
where $u_0:\R^+\to\mathbb{R}$ is a given function and $k=1,2.$ When $k\leq 3$, (2.1) is called the regular KdV generalized equation  while $k=4$ corresponds to the critical KdV equation, \cite{ larluc1,larluc2,martel, merle}. The regular case $k=1$ has been intensively studied  \cite{familark,larluc2}, while published results on initial-boundary value problems in generalized case $k\geq 2$ in unbounded domains, such as half-line, are almost unknown.

Hereafter subscripts $u_x,$ etc. denote the partial derivatives,
as well as $\partial_x$  when it is convenient.
By $(\cdot,\cdot)$ and $\|\cdot\|$ we denote the inner product and the norm in $L^2(\R^+),$
and $\|\cdot\|_{H^k(\R^+)}$ stands for the norm in $L^2$-based Sobolev spaces.

We will need the following result \cite{lady2}.
\begin{lemma}\label{lemma1}
Let $u\in H^1(\R^+)$ and $u(0,t)=0,$  then
\begin{equation}\label{2.5}
\|u\|_{L^4(\R^+)}\leq 2^{1/2}\| u_x\|^{1/2}\|u\|^{1/2},\;\;\|u\|_{L^8(\R^+)}\leq 4^{3/4}\| u_x\|^{3/4}\|u\|^{1/4}.
\end{equation}
\end{lemma}

\begin{proposition}\label{prop1}
	Let for a.e. fixed $t$  $u(x,t)\in H^1(\R^+).$ Then
	\begin{align}
		&\sup_{x\in{\R^+}}u^2(x,t)\le 2 \|u\|(t)\|u_x\|(t).
	\end{align}
\end{proposition}

\section{Existence theorem}\label{existence}

\begin{thm}\label{theorem1}
Given  $u_0\in H^3(\R^+)$  such that $u_0(0)=0.$ 
Then for all  positive $T$ there exists a unique regular solution to (2.1)-(2.3)
 such that
\begin{align*}
&u\in L^{\infty}(\R^+;H^3(\R^+))\cap L^2(\R^+;H^4(\R^+));\\
&u_t\in L^{\infty}(\R^+;L^2(\R^+))\cap L^2(\R^+;H^1(\R^+))
\end{align*}

\end{thm}

To prove this theorem, we will use regularization of (2.1)-(2.3) by an initial-boundary value problem on the half-line for the following Kawahara equation 

\begin{align}
	&L_{\epsilon}u^{\epsilon}\equiv u^{\epsilon}_t+u^{\epsilon k} u^{\epsilon}_x+u^{\epsilon}_{xxx}-\epsilon u^{\epsilon}_{xxxxx}=0\;\; \text{in}\;\; Q_T;\\
	&u^{\epsilon}(0,t) =\epsilon u^{\epsilon}_x(0,t)  =0,\; t>0;\\
	&u^{\epsilon}(x,0)=u^m_0(x),\ \ x\in\R^+,
\end{align}

where $\epsilon \in (0,1]$ and $u^m_0(x)$ is an approximation of $u_0(x)$ by functions $u^m_0\in H^5(\R^+)$ such that
$$\lim_{m\to \infty}\|u_0-u_0^m\|_{H^3(\R^+)}=0.$$
For $\epsilon=1$, (3.1)-(3.3) has been studied in \cite{famart} where it has been established that for $k< 8$  this problem 
has a regular solution for all $T>0$ provided a function $u^m_0(x)$ satisfies some compatibility conditions at $x=0.$ Problem (3.1)-(3.3) with $\epsilon \in (0,1)$ can be treated with the same arguments.\\
Our goal is to prove necessary a priori estimates for $u^{\epsilon}(x,t)$       independent of $\epsilon$ that will allow us to pass to the limit as $\epsilon\to 0$ and $m\to \infty$  getting a unique global regular solution for $k=1,2.$
\begin{proof}

{\bf Estimate I.}
Multiply (3.1) by $u^{\epsilon}$ and integrate over $\Omega \times (0,t)$ to obtain

\begin{align} &\| u^{\epsilon} \|^2(t) +\epsilon \int_0^t  |u^{\epsilon}_{xx}(0,\tau)|^2\, d\tau  = \|u^{\epsilon}_0 \|^2\leq \|u_0\|^2,\;\;t>0.
\end{align}

{\bf Estimate II.}\label{2-d estimate}
Write the inner product
$$2\left(L_{\epsilon}u^{\epsilon},(1+x)u^{\epsilon}\right)(t)=0,
$$
dropping the index $\epsilon$, in the form: 
\begin{align}
&\frac{d}{dt}\left((1+x),u^2\right)(t)
+{\epsilon}u_{xx}^2(0,t)
+3\|u_x\|^2(t)\notag\\&+5\epsilon\|u_{xx}\|^2(t)=-2((1+x)u^k u_x,u)(t).
\end{align}
Taking into account (2.4) and (3.4), we obtain
\begin{align*}
I=&-2((1+x)u^k u_x,u)(t)=\frac{4}{k+2}(u^ku^2)(t)\\&\leq 2sup_{(\R^+)}u^2(x,t)\|u\|^2(t)\leq 4\|u\|^3(t)\|u_x\|(t)\\&\leq
2[\|u_x\|^2(t)+\|u_0\|^6].
\end{align*}

Substitutinjg $I$ into (3.5), we get

\begin{align*}
&\frac{d}{dt}\left((1+x),u^2\right)(t)
+{\epsilon}u_{xx}^2(0,t)
+\|u_x\|^2(t)\notag\\&+5\epsilon\|u_{xx}\|^2(t)\leq \|u_0\|^6.
\end{align*}

Hence

\begin{align}
&\left((1+x),u^2\right)(t)
+\int_0^t\big[{\epsilon}u_{xx}^2(0,\tau)
+\|u_x\|^2(\tau)\notag\\&+5\epsilon\|u_{xx}\|^2(\tau)\big]d\tau\leq T\|u_0\|^6 +((1+x),u^2_0) \notag\\&
\leq(1+T\|u_0\|^4)((1+x),u^2_0)\equiv C_1((1+x),u^2_0).
\end{align}

{\bf Estimate III.}\label{3-d estimate}
Write the inner product
$$2\left(L_{\epsilon}u^{\epsilon},(1+x)^2u^{\epsilon}\right)(t)=0,
$$
dropping the index $\epsilon$, in the form: 
\begin{align}
&\frac{d}{dt}\left((1+x)^2,u^2\right)(t)
+2{\epsilon}u_{xx}^2(0,t)
+6\|(1+x)^{1/2}u_x\|^2(t)\notag\\&+10\epsilon\|(1+x)^{1/2}u_{xx}\|^2(t)=-2((1+x)^2u^k u_x,u)(t)
\end{align}
Taking into account (2.4) and (3.4), we obtain
\begin{align*}
I=&-2((1+x)^2u^k u_x,u)(t)=-\frac{4}{k+2} ((1+x)u^ku^2)(t)\\&\leq\frac{4}{k+2}sup_{(\R^+)}u^2(x,t)((1+x),u^2)(t)\leq \frac{8}{3}\|u\|(t)\|u_x\|(t)((1+x),u^2)(t)\\&\leq
\frac{4}{3}\|(1+x)^{1/2}u_x\|^2(t)+\frac{4}{3}\|u_0\|^2((1+x),u^2)^2(t).
\end{align*}

Substituting $I$ into (3.5), we obtain

\begin{align*}
	&\frac{d}{dt}\left((1+x)^2,u^2\right)(t)
	+2\|(1+x)^{1/2}u_x\|^2(t)\notag\\&+10\epsilon\|(1+x)^{1/2}u_{xx}\|^2(t)\leq 2\|u_0\|^2((1+x)^2,u^2)^2(t)\notag\\&
	\leq 2C^2_1\|u_0\|^2((1+x),u^2_0)^2.
\end{align*}

This implies
\begin{align}
&((1+x)^2,|u^{\epsilon}|^2)(t)\notag\\
&+\int_0^t\Big[2\|(1+x)^{1/2}u^{\epsilon}_x\|^2(\tau)+10\epsilon\|(1+x)^{1/2}u^{\epsilon}_{xx}\|^2(\tau)\Big]\,d\tau\notag\\&\leq \big(1+2C^2_1T\|u_0\|^2((1+x)^2,u_0^2)^1\big]((1+x)^2,u^2_0)\notag\\&\equiv C_2((1+x)^2,u^2_0).
\end{align}

{\bf Estimate IV}\\
Dropping the index $\epsilon$, write the inner product
$$
2\left((1+x)u^{\epsilon}_t),\partial_t(L_{\epsilon}u^{\epsilon}\right)(t)=0
$$
as
\begin{align}\label{3.13}
	\frac{d}{dt}
	&\left((1+x),u_t^2\right)(t)+\epsilon u_{xxt}^2(0,t)+3\|u_{xt}\|^2(t)	+5\epsilon\|u_{xxt}\|^2(t)\notag\\
	&=2\left((1+x)u^ku_t,u_{xt}\right)(t)+2(u^k,u_t^2)(t).
\end{align}
Taking $k=2$, we estimate
\begin{align*}
	I_1
	&=2\left((1+x)u^ku_t,u_{xt}\right)(t)\\&\leq 2\|u_{xt}\|(t)\sup_{(\R^+)}|u(1+x)^{1/2}u(x,t)|^2\|(1+x)^{1/2}u_t\|(t)\\
	&\le 4\|u_{xt}\|(t)\|(1+x)^{1/2}u\|(t)\|((1+x)^{1/2}u)_x \|(t)\|(1+x)^{1/2}u_t\|(t)\\
	&\le 2\| u_{xt}\|^2(t)+2\|(1+x)^{1/2}u\|^2(t)\|\frac{u}{2(1+x)^{1/2}}\\&+(1+x)^{1/2}u_x\|^2(t)((1+x)^{1/2},u^2_t)(t)\\ &
	\le 2\| u_{xt}\|^2(t)+4\|(1+x)^{1/2}u\|^2(t)\Big[\|u\|^2(t)\\&+((1+x),u^2_x)(t)\Big]((1+x),u^2_t)(t).		
\end{align*}

Similarly,
\begin{align*}
	I_2
	&=2(u^k,u_t^2)(t)\le 2\sup_{\R^+}u^2(x,t)((1+x),u^2_t)(t)\\
	&\leq 4\|u\|(t)\|u_x\|(t)((1+x),u^2_t)(t)\\
	&\leq 2\Big(\|u_0\|^2+\|u_x\|^2(t)\Big)((1+x),u^2_t)(t).
\end{align*}
Substituting $I_1, I_2$ into (3.9), we come to the inequality
\begin{align}\label{3.14}
&\frac{d}{dt}
((1+x),u_t^2)(t) +\epsilon u^2_{xxt}(0,t)+\|u_{xt}\|^2(t)	+5\epsilon \|\ u_{xxt}\|^2(t)\notag\\& 
 \leq C\big[\|(1+x)^{1/2}u_0\|^2+\|(1+x)^{1/2}u_x\|^2(t)\Big]((1+x),u^2_t)(t).
\end{align}

By (3.8), $\|(1+x)^{1/2}u_x\|^2(t) \in L^1(0,T),$ hence

\begin{equation} 
\|u_t\|^2(t)\leq ((1+x),u^2_t)(t)\leq C((1+x),|u^{\epsilon}_t|^2)(0),
\end{equation}

where $C$ depends on $\|u_0\|,T$ and
\begin{align}
	&((1+x),|u^{\epsilon}_t|^2)(0)=\int_{\R^+} (1+x)\Big[-u_{0xxx}-u^k_0 u_{0x}+\epsilon\partial^5_x u_0\Big]^2 dx\notag\\&
\leq 2\int_{\R^+} (1+x)\Big[-u_{0xxx}-u^k_0 u_{0x}\Big]^2 dx +2\epsilon^2\int_{\R^+}(1+x)\Big[ \partial^5_x u_0\Big]^2dx.	
\end{align}

Returning to (3.10), we find

\begin{align}
	&((1+x),|u^{\epsilon}_t|^2)(t)+\int_0^T\{\|u^{\epsilon}_{xt}\|^2(t)+\epsilon\|u^{\epsilon}_{xxt}\|^2(t)\}dt\notag\\&	\leq	C((1+x),|u^{\epsilon}_t|^2)(0).
\end{align}

{\bf Passage to the limit as $\epsilon \to 0.$}

Rewrite (3.1) in the form

\begin{align}& (u^{\epsilon}_t,\phi)(t)+(u^{\epsilon}_x,\phi_{xx}(t)) -\frac{1}{k+1}(u^{\epsilon (k+1)},\phi_x)(t)\notag\\&
	+\epsilon(u^{\epsilon}_{xx},\phi_{xxx} )(t)=0,
\end{align}   

where $\phi(x)$ is an arbitrary function such that $$\phi\in H^3(\R^+);\;\;\phi(0)=\phi_x(0)=\phi_{xx}(0)=0.$$
Making use of (3.8) and (3.13), we get $$ \epsilon^{1/2}\|u^{\epsilon}_{xx}\|(t)\leq C,\;\;\|u^{\epsilon}_x\|(t)\leq C,$$ hence
$$ \lim_{\epsilon\to 0} \epsilon\|u^{\epsilon}_{xx}\|(t)\|\phi_{xxx}\|=0,$$ and we obtain from (3.14)

\begin{equation}
 (u_t,\phi)(t)+(u_x,\phi_{xx})(t) -\frac{1}{k+1}(u^{(k+1)},\phi_x)(t)=0.
\end{equation}

Moreover,
\begin{align}&\lim_{\epsilon\to 0}((1+x),|u^{\epsilon}_t|^2)(0)=((1+x),|u^m_t|^2)(0)\notag\\&
=\int_{\R^+}(1+x)\{-u^m_{0xxx}-u^{mk}u^m_{0x}\}^2 dx,
\end{align}
consequently,
$$ \lim_{m\to\infty}\int_{\R^+}(1+x)\{-u^m_{0xxx}-u^{mk}u^m_{0x}\}^2 dx=\int_{\R^+}(1+x)\{-u_{0xxx}-u^{k}u_{0x}\}^2 dx.$$

Multiplying (3.1) by $(1+x)u$ and taking into account (3.7),(3.13) we come to the inequlaity
$$ \|u_x\|^2(t)\leq 2\|u_0\|^6+2|((1+x)u_t,u)(t)|\leq C,$$
then 

$$ u_{xxx}=-u^ku_x-u_t\in L^{\infty}((0,T);L^2(\R^+)).$$
In turn, multiplying this equality by $ -2((1+x)u_{xx},$ we find that\\ $u_{xx}\in L^{\infty}((0,T);L^2(\R^+)),$ hence 
$$ u\in L^{\infty}((0,T);H^3(\R^+)) ,\;\; u_t\in  L^{\infty}((0,T);L^2(\R^+))\cap  L^2((0,T);H^1(\R^+)).$$
Since $ u^ku_x\in  L^2((0,T);H^1(\R^+)),$ then
$ u_{xxx}\in  L^2((0,T);H^1(\R^+))$ and
\begin{align*}
	&u\in  L^{\infty}((0,T);H^3(\R^+))\cap  L^2((0,T);H^4(\R^+));\\&
	u_t\in  L^{\infty}((0,T);L^2(\R^+))\cap  L^{\infty}((0,T);H^1(\R^+)).
	\end{align*}

This proves the existence part of Theorem 3.1.

{\bf Uniqueness}
Let $u_1, u_2$ be two distinct solutions to (3.1)-(3.3) and $z=u_1-u_2.$ Considering more interesting case $k=2$, we find
$$Lz=z_t+z_{xxx}=-\frac{1}{3}(u^3_1-u^3_2)_x=-\frac{1}{3}\{z(u^2_1+u_1u_2+u^2_2)\}.$$

Multiplying this equation by $2(1+x)z$, we obtain
\begin{align}&\frac{d}{dt}((1+x,z^2)(t)+3\|z_x\|^2(t)+z^2_x(0,t)\notag\\&=-\frac{2}{3}(z[u^2_1+u_1u_2+u^2_2])_x,(1+x)z)(t).
	\end{align}

Estimate
\begin{align*}&I=-\-\frac{2}{3}(z[u^2_1+u_1u_2+u^2_2])_x,(1+x)z)(t)=\frac{1}{3}((u^2_1+u_1u_2+u^2_2),z^2)(t)\\&
+	\frac{1}{3}([(2u_1+u_2)u_{1x}+(2u_2+u_1)u_{2x}],(1+x)z^2)(t)\\&
\leq	\frac{2}{3}((u^2_1+u^2_2),z^2)(t)+\frac{1}{6}([(2u_1+u_2)^2+(2u_2+u_1)^2\\&+u^2_{1x}+u^2_{2x}],(1+x)z^2)(t)
	\leq \frac{2}{3}[\sup_{\R^+}u^2_1+\sup_{\R^+}u^2_2]((1+x),z^2)(t)\\&
	+ \frac{4}{3}\sup_{\R^+}[u^2_1+u^2_2+u^2_{1x}+u^2_{2x}]((1+x)z^2)(t)\\&
	\leq C\{\|u_1\|^2_{H^2(\R^+)}(t)+\|u_2\|^2_{H^2(\R^+)}(t)\}((1+x),z^2)(t)\leq C((1+x),z^2)(t).
	\end{align*}
	
Substituting $I$ into (3.17), we come to the inequality
$$\frac{d}{dt}((1+x),z^2)(t)\leq C((1+x),z^2)(t).$$
Since $z(x,0)\equiv 0,$,then
$$ \|z\|^2(t)\leq ((1+x),z^2)(t)\equiv 0\;\; t>0.$$
This completes the proof of Theorem 3.1
\end{proof}
\begin{remark}In the regular case $k=3$, it is possible to prove local in $t$ the existence and uniqueness of regular solutions. On the other hand, there are papers, \cite{lipaz}, where a local dissipation term $au,\;a>0$ is added to (1.1). It helps to prove the existence, uniqueness and exponential decay of small solutions.
	\end{remark}

{\bf Conclusions.} An initial-boundary value problem for the  generalized KdV equation posed on the right half-line  has been considered.The existence and uniqueness of a regular global solution for all positive $T$ and arbitrary smooth initial data have been established.

\end{document}